# Una tentazione affascinante

## Sull'infinito in matematica

di Claudio Bernardi

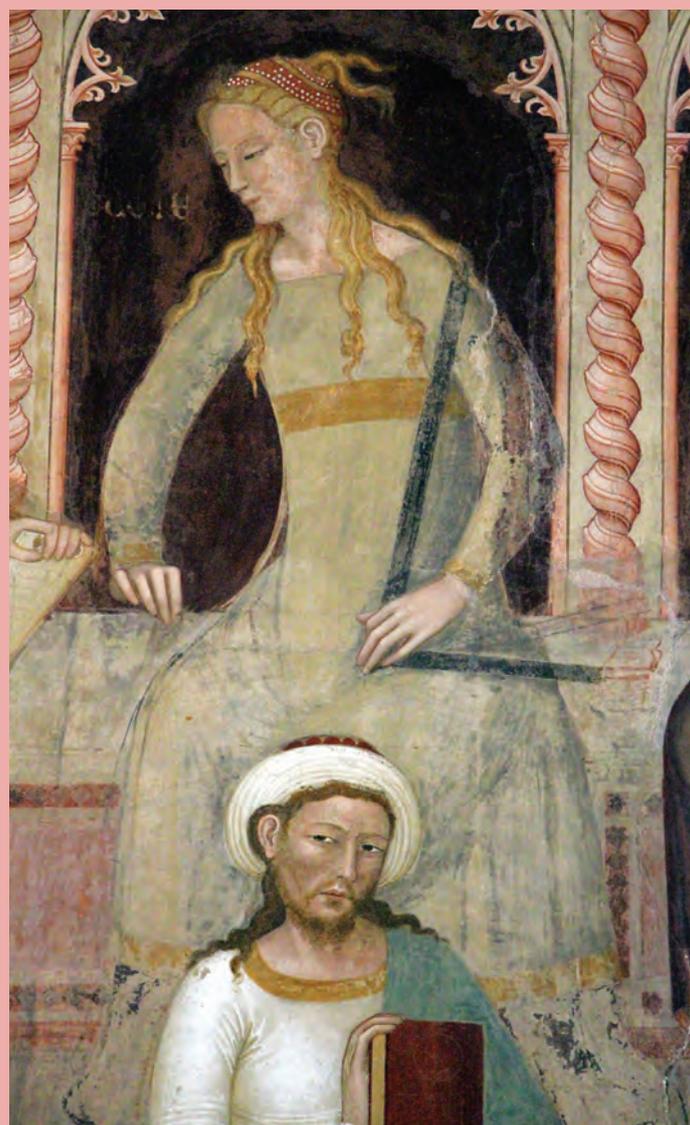

a.
Affresco di Andrea di Bonaiuto, raffigurante Euclide, nel cappellone degli Spagnoli, l'antica sala capitolare della chiesa di Santa Maria Novella a Firenze. Euclide si trova ai piedi della "Geometria". L'affresco fa parte del "Trionfo di San Tommaso d'Aquino": nella parte inferiore si trovano quattordici stalli decorati, nei quali siedono le personificazioni delle sacre scienze e delle arti liberali, ai piedi di ciascuna delle quali si trova un illustre rappresentante.

Secondo Hermann Weyl, matematico tedesco della prima metà del '900, "la matematica è la scienza dell'infinito". Su questa affermazione si può discutere, ma sta di fatto che, a differenza di quanto avviene in altri contesti scientifici, in matematica si parla con tranquillità di infinito.
Fin dalle prime classi della scuola primaria si insegna che ad ogni numero naturale si può aggiungere 1, in modo da ottenere un numero più grande di quello iniziale. Già a questo punto, il bambino capisce che ci sono infiniti numeri naturali. Quando poi si comincia la geometria, si spiega che una retta assomiglia a un filo teso, ma va pensata infinitamente lunga. Un po' più difficile è far capire agli studenti la densità: tra due punti (così come fra due numeri razionali o reali) ne è compreso almeno un altro, e quindi infiniti altri.
Ma il rapporto dei matematici con l'infinito non è stato sempre tranquillo. Basta ricordare che per Euclide "retta" significava "segmento prolungabile": in altre parole Euclide accettava un infinito "potenziale" (un segmento si può allungare quanto si vuole, ma ogni volta abbiamo a che fare con un segmento), ma non l'infinità di una retta in senso attuale. Le difficoltà con l'infinito nascevano dai paradossi a cui l'accettazione dell'esistenza di insiemi infiniti sembrava inevitabilmente condurre.
Galileo rimane sconcertato dal fatto che i quadrati perfetti sembrano tanti quanti i numeri naturali (ad ogni naturale corrisponde un quadrato, e viceversa), mentre, d'altro lato, i quadrati sono solo una piccola parte dell'insieme dei naturali: anzi, via via che crescono i valori che si considerano, i quadrati diventano sempre più "radi" nell'insieme dei numeri naturali. Galileo si esprime nel modo seguente, quasi suggerendo che è meglio rinunciare a uno studio matematico dell'infinito: "Queste son di quelle difficoltà che derivano dal discorrer che noi facciamo col nostro intelletto finito intorno a gli infiniti, dandogli quegli attributi che noi diamo alle cose finite e terminate; il che penso che sia inconveniente, perché stimo che questi attributi di maggioranza, minorità ed ugualità non convenghino a gl'infiniti, de i quali non si può dire, uno essere maggiore o minore o eguale all'altro".

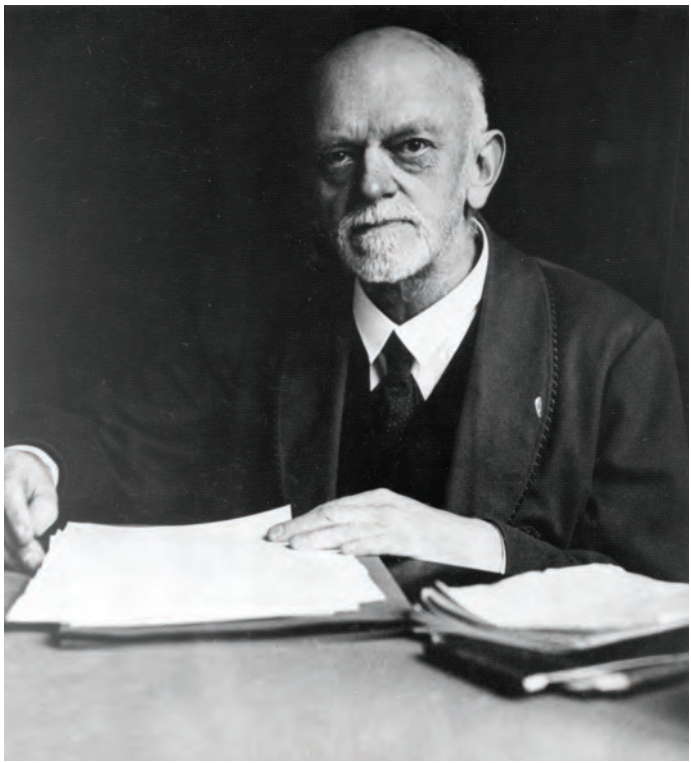

Anche il grande matematico Carl Friedrich Gauss ebbe a dire: "Io devo protestare nel modo più deciso contro l'uso dell'infinito come qualcosa di compiuto, cosa che non è permessa in matematica. L'infinito non è che una *façon de parler* (modo di dire, ndr)".
Ma poi, verso la fine dell'800 arrivò Georg Cantor: nella teoria degli insiemi di Cantor si riescono a studiare e confrontare insiemi infiniti, senza cadere in paradossi. Questo permette, oggi, di trattare senza paura insiemi infiniti di numeri (come quello dei numeri reali *R* e quello dei numeri complessi *C*), di punti, di funzioni, ecc.

Sempre nella seconda metà dell'800, a opera di Karl Weierstrass e altri, si sistemano in modo rigoroso i concetti alla base dell'analisi matematica: con le celebri formule del tipo "$\forall \varepsilon > 0 \, \exists \delta > 0 \ldots$" (che si legge "per ogni ε maggiore di 0 esiste un δ maggiore di 0", ndr) si danno definizioni chiare per scritture come $\lim_{x \to \infty} f(x) = L$ (ovvero "limite di $f(x)$ per $x$ che tende a infinito è uguale a $L$", ndr) e $\lim_{x \to \infty} f(x) = \infty$ ("limite di $f(x)$ per $x$ che tende a infinito è uguale a infinito", ndr), e si definisce la "continuità" di una funzione, rinunciando a discutibili variazioni infinitesime delle variabili.
Oggi, in matematica si parla con tranquillità di infinito. Il ricorso all'infinito rende spesso addirittura più semplice lo studio della matematica: in particolare, la teoria delle derivate e degli integrali è più semplice e più potente della teoria delle differenze finite, che si riferisce a funzioni che hanno per dominio un insieme discreto invece che tutto l'insieme *R* dei numeri reali.
In altri casi, accettando l'infinito si ottengono strutture più regolari e più eleganti. Così, aggiungendo a un piano elementare i punti all'infinito (punti "impropri") si ottiene un piano "proiettivo", dove si enunciano proprietà generali senza eccezioni e senza che sia necessario distinguere vari casi: due rette distinte hanno sempre uno e un solo punto in comune, due coniche "non degeneri" si possono sempre trasformare una nell'altra, ecc.
Vale la pena di osservare un altro legame fra matematica e infinito, in parte condiviso da altre scienze: una dimostrazione sostituisce infinite verifiche. Pensiamo per esempio al teorema sulla somma degli angoli interni di un triangolo. Per verificare che la somma è uguale a un angolo piatto in tutti i triangoli, anche se avessimo a disposizione strumenti di misura perfetti, dovremmo eseguire infiniti controlli. Invece, una singola dimostrazione ci convince che il teorema è vero in generale, rendendo superflui gli infiniti controlli.
Come disse il grande matematico David Hilbert nel 1921, "L'infinito! Nessun altro problema ha mai scosso così profondamente lo spirito umano; nessuna altra idea ha stimolato così proficuamente il suo intelletto; e tuttavia nessun altro concetto ha maggior bisogno di chiarificazione che quello di infinito".

b.
Un ritratto del 1937 del grande matematico David Hilbert.

**Biografia**
**Claudio Bernardi** è professore ordinario di matematica alla Sapienza Università di Roma. Si occupa di logica matematica e di fondamenti della matematica e anche di didattica della matematica. È stato docente alla Ssis e al Tfa, è il responsabile del piano lauree scientifiche per la matematica alla Sapienza e collabora con l'Accademia dei Lincei per iniziative legate all'insegnamento.

**Link sul web**

www.dm.unipi.it/~berardu/Didattica/Appunti/scuole060210.pdf

www.treccani.it/vocabolario/infinito/

http://ed.ted.com/lessons/how-big-is-infinity

http://maddmaths.simai.eu/divulgazione/matematica-indispensabile/
il-concetto-matematico-di-cui-non-potremmo-fare-a-meno-infinito-di-sandra-lucente/